\documentclass[12pt]{article}

\usepackage{amsmath}
\usepackage{multirow} 
\usepackage[cp850]{inputenc}
\usepackage{latexsym,amsfonts}
\usepackage{graphicx}
\usepackage{subcaption}
\usepackage{pdflscape}
\usepackage{caption}
\usepackage{tikz-cd} 
\usepackage{tikz}
\usepackage{enumitem} 
\usepackage{authblk}
\newtheorem{thm}{Theorem}[section]
\newtheorem{note}{Note}

\newtheorem{cor}[thm]{Corollary}

\date{}
\author{Vladimir D. Tonchev}
\affil{Department of Mathematical Sciences, Michigan Technological 
University, Houghton, MI USA 49931, tonchev@mtu.edu}

\title{ On Pless symmetry codes, ternary QR codes,  and related Hadamard matrices and designs}

\begin{document}
\maketitle

\begin{abstract}

It is proved that a code $L(q)$ which is monomially
 equivalent to  the Pless symmetry
code  $C(q)$ of length $2q+2$ contains the (0,1)-incidence matrix of a Hadamard 3-$(2q+2,q+1,(q-1)/2)$  design 
$D(q)$ associated with
a  Paley-Hadamard matrix of type II. Similarly,
any ternary extended quadratic residue code contains
the incidence matrix of a Hadamard 3-design  associated with a
Paley-Hadamard matrix of type I. If $q=5, 11, 17, 23$,
then the full permutation automorphism group of $L(q)$
coincides with the full automorphism group of $D(q)$, and a similar result holds for the ternary extended quadratic residue
codes of lengths 24 and 48.
All Hadamard matrices of order 36 formed by codewords
of the Pless symmetry code $C(17)$  are 
enumerated and classified up to equivalence.
There are two equivalence classes of such matrices:
the Paley-Hadamard matrix  $H$ of type I
with a full automorphism group of order 19584, and a second
regular Hadamard matrix $H'$ such that the symmetric
2-$(36,15,6)$ design $D$ associated with $H'$ has trivial
full automorphism group, and the incidence matrix of $D$ spans
a ternary code equivalent to $C(17)$.

{\bf Keywords:} Pless symmetry code, Hadamard matrix,
Hadamard 3-design, Hadamard 2-design,
Paley-Hadamard matrix.

\end{abstract}

\section{Introduction}

We assume familiarity with the basic facts and notions 
from error-correcting codes and combinatorial designs and Hadamard matrices \cite{AKb}, \cite{BJL}, \cite{Hall}, \cite{HP}.
All codes in this paper are ternary.
A monomial matrix with entries from $GF(3)$ is a square matrix
such that every row and every column contains exactly one nonzero entry.
 By  an automorphism group of a
ternary  code we mean {\it monomial} automorphism group, unless specified otherwise. The {\it permutation} automorphism group of a code
is the subgroup of  its monomial automorphism group that
consists of coordinate permutations only.

A Hadamard matrix of order $n$ is an $n \times n$ matrix $H$ of $1$'s and $-1$'s such that $HH^T=nI$, where $I$ is the identity matrix. 
It follows that $n=1, 2$, or $n=4t$ for some integer $t\ge 1$.
An automorphism of a Hadamard matrix $H$ is a pair of 
$\{ 0,1,-1\}$-monomial matrices $L$, $R$ such that $LHR=H$.
Two Hadamard matrices $H_1$, $H_2$ of the same order are {\it equivalent}
if there are monomial matrices $L$, $R$ such that $LH_{1}R=H_2$.
A Hadamard matrix $H$ is {\it normalized} with respect to its $i$th row and $j$th 
column if all entries in row $i$  and column $j$ are equal to 1. 
If $H$ is a Hadamard matrix of order $n=4t$ that is normalized with respect to row
 $i$ and column $j$, deleting the $i$th row and the $j$th column and replacing all $-1$'s with zeros gives the $(0,1)$-incidence matrix of a symmetric 2-$(4t-1,2t-1,t-1)$ design $D$
  called a Hadamard 2-design, while deleting the $j$th column of $H$ and the $j$th column of 
  $-H$ from the matrix $(H,-H)$ gives the point-by-block $(\pm 1)$-incidence matrix of a 
  3-$(4t,2t,t-1)$ design $D^*$, called a Hadamard 3-design obtained from $H$ with respect to column $j$. The design $D$ is the {\it derived} design of $D^*$ with respect to its $i$th point.
A Hadamard matrix $H$ of order $n=4t$ is {\it regular} if all rows of $H$ contain the same number $k$ of $-1$'s. It follows that $t=m^2$ for some integer $m$, $k=2m^2 \pm m$, and 
replacing all $-1$'s with zeros gives the $(0,1)$-incidence matrix of a symmetric 
2-$(4m^2,2m^2\pm t,m^2\pm m)$ design.
For more on Hadamard matrices and related designs, see, for example,
 \cite[Chapter 7]{AKb}, \cite[Chapter 14]{Hall}, \cite[Sec. 8.9]{HP}.

Let $q$ be an odd prime power such that 
$q \equiv -1 \pmod 3$. The Pless symmetry code $C(q)$
\cite{Pless69}, \cite{Pless72} of length $n=2q+2$ is a ternary 
self-dual code with a generator matrix 
 
 \begin{equation}
 \label{eq1}
 G=(I_{q+1}, S_q), 
 \end{equation}
 where $I_{q+1}$ is the identity matrix of order $q+1$, and 
$S_q =(s_{i,j})$ is a 
$(q+1)\times (q+1)$  matrix defined as follows. The rows and columns of $S_q$ are labeled by $\infty$ and the $q$ elements of the finite field $GF(q)$ of order $q$, 
where $s_{\infty,\infty}=0$,
$s_{a,a}=0$,  $s_{\infty,a}=1$ for $a\in GF(q)$,
 $s_{a,\infty}=1$ if $-1$ is a square in $GF(q)$,
 and $s_{a,\infty}=-1$ if $-1$ is not a square in $GF(q)$ for 
$a\in GF(q)$, and $s_{a,b} =1$ for  $a, b\in GF(q)$ such
that $a\neq b$ and $a-b$ is a square in $GF(q)$,
and $s_{a,b} =-1$ for  $a, b\in GF(q)$ such
that $a\neq b$ and $a-b$ is not a square in $GF(q)$.
For example,  if $q=5$, the rows and columns of $S_5$ are labeled by
$\infty, 0, 1, 2, 3, 4$, and
\[ S_5 =\left(\begin{array}{rrrrrr}
0 & 1 & 1 & 1      & 1  & 1 \\
1 & 0 & 1 & -1 & - 1  & 1\\
1 &  1  & 0 & 1  & -1 & -1\\
1 & -1 & 1 & 0 & 1 & -1\\
1 & -1 & -1 & 1 & 0 & 1\\
1 & 1 & -1 & -1 & 1 & 0
\end{array}
\right).
\]
The main property of $S_q$ is that $S_{q}S_{q}^T=qI_{q+1}$,
which implies  $S_{q}S_{q}^T \equiv -I_{q+1} \pmod 3$;
hence $C(q)$ is self-dual.
The  symmetry codes
$C(5)$\footnote{The symmetry code for $q=5$ is equivalent to the extended
 ternary Golay code.} ,
$C(11)$, $C(17)$, $C(23)$ and $C(29)$
are extremal self-dual codes of length $n$ divisible by 12
and minimum distance $d$ meeting the
Mallows-Sloane  upper bound 
$d\le 3[n/12]+3$ \cite{MS};
hence these codes support 5-designs by the Assmus-Mattson
Theorem \cite{AM}.

Pless \cite{Pless72} proved that in addition to the
trivial monomial automorphism of order 2 corresponding to 
the negation of all code coordinates, 
the monomial automorphism group of $C(q)$
contains a subgroup of order $q(q^2 -1)$
isomorphic to $PGL_{2}(q)$.
In addition,
it was proved in  \cite[Theorem 4.2]{Pless72}   that
the symmetry code $C(q)$ contains
a set of $2q+2$ codewords  of weight $2q+2$
that form a Hadamard matrix  of order $2q+2$.
Specifically, if $q\equiv 1 \pmod 4$ the Hadamard matrix
formed by codewords of weight $2q+2$  is
\begin{equation}
\label{eq2}
 H_{1}(q) =\left(\begin{array}{rr}
 I + S_q & -I +S_q\\
 -I + S_q & -I - S_q
 \end{array}\right),
 \end{equation}
while if $q\equiv -3 \pmod 4$ the Hadamard matrix is
\begin{equation}
\label{eq3}
H_3(q) = \left(\begin{array}{rr}
 I + S_q & I +S_q\\
  I - S_q & -I + S_q
 \end{array}\right),
  \end{equation}
  where $I$ is the identity matrix of order $q+1$.

The Hadamard matrices (\ref{eq2}), (\ref{eq3})
are known in the combinatorial literature as Paley-Hadamard matrices of type II
 \cite[14.1]{Hall}, \cite{LS}, \cite{Paley}.
If $q\equiv 3 \pmod 4$,  the $(q+1)\times(q+1)$ matrix
obtained by bordering the matrix $S_q -I$ with one all-one
row and one all-one column is a Hadamard matrix of order $q+1$,
known as a Paley-Hadamard matrix, or a Paley-Hadamard matrix
of type I \cite[14.1]{Hall}, \cite{LS}, \cite{Paley}.
The unique (up to equivalence) Hadamard matrix
of order 12 is both a Paley-Hadamard matrix of type I for $q=11$
and a Paley-Hadamard matrix of type II for $q=5$, and its full automorphism group modulo its center of order 2 is  the Mathieu group $M_{12}$ (Hall \cite{Hall62}).
The full automorphism group of a Paley-Hadamard matrix
of type I  for  $q>11$ was determined by
Kantor \cite{K} and is of order $q(q^2-1)$, while
the full automorphism group of a Paley-Hadamard matrix of type 
II for $q>5$ was determined by de Launey and Stafford \cite{LS}, 
and is of order $4fq(q^2 -1)$ if $q=p^f$, where $p$ is prime. 

If $q=5$, 11, or 23, the number
of all codewords of full weight $2q+2$ in the symmetry code $C(q)$
is exactly $4q+4$ \cite{Pless72}.
These codewords span the code and consist of the rows of the  
Hadamard matrix $H_1(q)$ from (\ref{eq2}) 
(resp. $H_3(q)$ from (\ref{eq3})) and  its negative, or $2H_1(q)$  
(resp. $2H_2(q)$); hence  the full monomial automorphism group of $C(q)$ coincides  
with the full automorphism group of $H_1(q)$ (resp. $H_3(q)$) \cite{Pless72}.

The symmetry code $C(17)$ of length 36 contains
888 codewords of weight 36; hence it is not clear
whether the automorphism group of the Hadamard
matrix   $H_1(17)$ (\ref{eq2}) is 
 the full automorphism group of $C(17)$.
In  Section 2, we prove that the full automorphism group
of $C(17)$ coincides with the full automorphism group
of $H_1(17)$, the latter being
a Paley-Hadamard matrix of type II; hence its  order
is  $4\cdot{17}(17^2 -1)=19584$.
In addition, we classify all Hadamard matrices of order 36
having as rows codewords of $C(17)$ of weight 36, and show that
up to equivalence, there are exactly two such matrices:
 $H_1(17)$,
and a second Hadamard  matrix $H'$ having the property that all
Hadamard 3-$(36,18,8)$ designs associated with $H'$ are
isomorphic and have trivial full automorphism group 
of order 1.
The full automorphism group of $H'$ is of order $72$ and 
is transitive on the set of 72 rows (as well as the set of 72 columns)
of $H'$ and $-H'$. The 3-rank of $H'$ is 18; thus $C(17)$
is the row space of $H'$. The Hadamard matrix
$H'$ is regular, and the symmetric 2-$(36,15,6)$ design
with $\pm 1$-incidence matrix $H'$ has trivial full automorphism
group.

In Section \ref{s3} we discuss  Paley-Hadamard matrices of Type I and 
Hadamard 3-designs arising from  extended ternary quadratic residue
codes.

\section{Hadamard matrices  and designs arising from symmetry codes}
\label{s2}
The sum (over $GF(3)$) of all rows of the generator matrix (\ref{eq1}) of 
the symmetry code $C(q)$ is a vector $v$
of full Hamming weight $2q+2$, with all components equal to 1
if $-1$ is  not a square in $GF(q)$, (that is, $v$ is the constant
all-one vector $\bar{1}=(1,\ldots,1)$), and $v$ has $2q+1$ components
equal to 1, and the position labeled by $\infty$ is equal
to $-1$ whenever $-1$ is a square in $GF(q)$.

Next, we  consider a code which is monomially equivalent to the Pless symmetry code $C(q)$,
and always contains the all-one vector, namely the code $L(q)$
with a generator matrix $G'$  given by 
\begin{equation}
\label{eq4}
 G' =(I_{q+1}, U_q), 
 \end{equation}
 where $U_q$ is a $(q+1)\times(q+1)$ matrix
obtained from $S_q$ by replacing
every nonzero entry in the column labeled by $\infty$
with $-1$. 
Clearly, the generator matrix $G'$ (\ref{eq4}) is identical with
the generator matrix $G$ (\ref{eq1})  if $-1$ is not  a square in $GF(q)$, and is obtained
by negating one column of $G$ if $-1$ is a square in $GF(q)$.
Thus, the matrices (\ref{eq1}) and (\ref{eq4}) generate
monomially equivalent ternary codes.

Using (\ref{eq4}), we obtain a  parity check $P$ matrix
of $L(q)$ given by  
\begin{equation}
\label{eq5}
P = (-U_q^T, I_{q+1}).
\end{equation}
Note that since $L(q)$ is self-dual, the rows of $P$ are codewords 
of  $L(q)$. 
It is easy to check that the matrix $H$ given by
\begin{equation}
\label{eq6}
H=\left(\begin{array}{c}
G'+P\\
G' -P
\end{array}\right) = \left(\begin{array}{cc}
I_{q+1}-U_q^T & U_q +I_{q+1}\\
I_{q+1}+U_q^T & U_q -I_{q+1}
\end{array}\right)
\end{equation}
is a Hadamard matrix of order $2q+2$,
with rows being codewords of $L(q)$.
\begin{thm}
\label{t1}
The code $L(q)$ contains a set of $4q+2$ (0,1)-codewords
of weight $q+1$ that form the block-by-point incidence matrix
of a Hadamard 3-$(2q+2,q+1,(q-1)/2)$ design $D(q)$
associated with a Paley-Hadamard matrix of type II.
\end{thm}
{\bf Proof}.
All entries in the first row of the Hadamard matrix $H$ (\ref{eq6})
are equal to 1; that is, $H$ is normalized with respect to its
first row, and consequently, all entries in the first row of $-H$
are equal to $-1$. Adding the constant codeword
$\bar{2}=(2,\ldots,2)$  with all
entries equal to 2 to every row of the matrix
\[ \left(\begin{array}{rr}
 H\\ -H
 \end{array}\right)
 \]
gives a $(0,1)$-matrix $M$ with all-zero first row, and all-one row
labeled by the first row of $-H$. Deleting the all-zero row and
the all-one row from $M$ gives a $(4q+2)\times (2q+2)$
$(0,1)$-matrix $A$, being the block-by-point incidence matrix
of a Hadamard 3-$(2q+2,q+1,(q-1)/2)$ design associated with the first row 
of $H$. 
Clearly, $H$ is equivalent to the
corresponding matrix (\ref{eq2}) or (\ref{eq3}); hence $H$ is
equivalent to a Paley-Hadamard matrix of type II.
$\Box$

\begin{thm}
\label{t2}
If $q=5, 11, 17, 23$, the code $L(q)$ contains exactly
$4q+2$ (0,1)-codewords of weight $q+1$, and every
such codeword is the incidence vector of a block of the
Hadamard 3-design $D(q)$ from Theorem \ref{t1}.
\end{thm}
{\bf Proof}.  Let $m$ denote the total number of (0,1)-codewords
of weight $q+1$ in $L(q)$. By Theorem \ref{t1}, $m\ge 4q+2$.
 If $v\in L(q)$ is a (0,1)-codeword of weight $q+1$
then $v+\bar{1}$ is a codeword of full weight $2q+2$,
having $q+1$ components equal to 1, and $q+1$ components
equal to 2. Adding the codewords $\bar{1}$ and 
$\bar{2} = 2\cdot\bar{1}$
gives $m+2 \ge 4q+4$ codewords of weight $2q+2$.
Since $C(q)$ contains exactly $4q+4$ codewords of weight
$2q+2$ for $q=5, 11, 23$ \cite{Pless72}, the statement is true in these cases. 

The case $q=17$ needs additional analysis because
the symmetry code $C(17)$ , as well as  its  equivalent code 
$L(17)$, contains 888 codewords of weight 36 \cite{Pless72}.
The set of  all codewords of weight 36 is easily computed with Magma
\cite{magma}.  This set comprises of the following codewords:
\begin{itemize}
\item
the 36 rows of the Hadamard matrix $H$ (\ref{eq6}),
one of the  rows being $\bar{1}$, and 35 rows with 18 
components equal to 1,
and 18 components equal to  $-1$ (note that $-1\equiv 2 \pmod 3$);
\item
the 36 rows of $2H$ that include $\bar{2}$  and 35 rows
 with 18  components equal to 1, and 18 components equal to 2;
\item
 a set $T$ of 408 codewords having 15 components equal to 1 and 21 components equal to 2;
 \item   a set
$2T$ of 408 codewords obtained by multiplying every codeword from $T$ 
by 2. 
\end{itemize}

Note that adding $\bar{2}$ to any $(0,1)$-codeword of weight 18
gives a codeword of weight 36 with 18 $1$'s and 18 $2$'s; hence
the code $L(17)$ contains exactly 70 $(0,1)$-codewords of weight 18
obtained by adding the codeword $\bar{2}$ to the rows of $H$ and $2H$,
and these 70  $(0,1)$-codewords form the incidence matrix of the
3-design $D(17)$ from Theorem \ref{t1}.
$\Box$
\begin{note}
{\rm
The code $L(29)$ contains 19606 (0,1)-codewords of weight 30. 
It is an  open question whether this set  contains the incidence matrices of any
Hadamard 3-$(60,30,14)$ designs that are not isomorphic to $D(29)$.
The number of codewords of weight 60 in $L(29)$ is 41184. It seems likely
that there may be Hadamard matrices of order 60 formed by codewords of
weight 60 that are not equivalent to the Paley-Hadamard matrix of type II.
}
\end{note}
\begin{cor}
\label{cor1}
If $q=5, 11, 17, 23$, the full permutation automorphism group of $L(q)$
coincides with the full automorphism group
of the Hadamard 3-design $D(q)$ from Theorem \ref{t1}.
\end{cor}
{\bf Proof}. The results of De Launey and Stafford \cite{LS}
and Norman \cite{N} imply that the full automorphism group of
$D(q)$ has order $q(q-1)$ if $q>5$ is  prime. The full automorphism
group of $D(5)$ has order 7920.

 Any derived design with respect to a point of a Hadamard
3-$(2q+2,q+1,(q-1)/2)$ design  $D$ is a symmetric Hadamard 2-$(2q+1,q,(q-1)/2)$
design $D'$ of order $q - (q-1)/2 =(q+1)/2$. Since $q\equiv -1 \pmod 3$, 3 divides
$(q+1)/2$. If 9 does not divide the order $(q+1)/2$  (which is true if $q=5$, 11, or 23),
the rank of the incidence matrix if $D'$ over $GF(3)$ (or the 3-rank of $D'$)
is equal to $q+1$ (see, for example, Assmus and Key \cite{AK},
\cite{AKb}), hence the 3-rank
of $D$ is $q+1$ and the code $L(q)$ is spanned by the incidence matrix of $D$.
If $q=17$, a direct computation shows that the 
3-rank  of $D(17)$  is 18, hence $D(17)$ spans the code $L(17)$.

$\Box$

\begin{thm}
\label{t3}
(i) The code $L(17)$ contains two equivalence classes of Hadamard matrices
of order 36 having as rows codewords of weight 36,
with representatives the
Hadamard matrix $H$ (\ref{eq6}), which is equivalent to a
Paley-Hadamard matrix of type II and has full automorphism group
of order $19584$, and a second Hadamard matrix $H'$,
being a regular Hadamard matrix such that
the  symmetric 2-$(36,15,6)$ design $D'$ with $(0,1)$-incidence matrix $(H' +J)/2$,
where $J$ is the $36\times 36$ all-one matrix,
 has a  trivial automorphism group.\\
(ii) The row span of the incidence matrix of the   2-$(36,15,6)$ design $D'$
 is an extremal
ternary $[36,18,12]$ code equivalent to the symmetry code $C(17)$.\\
(iii) The full automorphism group of the code $L(17)$ 
coincides with the full automorphism group $H$.
\end{thm}
{\bf Proof}.
(i) 
In the context of Hadamard matrices, we consider the
element 2 of $GF(3)$ as $-1$.
  Using the notation from the proof of Theorem \ref{t2},
we define a graph $\Gamma$ having as vertices the 408 codewords
from $T$, where two codewords 
$u, v \in T$  are adjacent in $\Gamma$ if and only
if the Hamming distance between $u$ and $v$ is 18, or equivalently,
the intersection of the supports of the (0,1)-vectors $\bar{2}-u$ and $\bar{2}-v$ is of size 6.
Replacing all entries equal to 2 by zero in every vector from $T$ gives
a set $T(0,1)$ of (0,1)-vectors of weight 15.
Using  the restricted Johnson bound, it is easy to verify that
the maximum number codewords in a binary constant weight code of length 36
with  codewords of weight 15 and  minimum distance 18, is 36.
Every set $K$ of  36 vectors from $T(0,1)$ that meets the Johnson bound
 corresponds to a clique of size 36 in $\Gamma$, and the $36\times 36$
matrix having as rows the vectors from $K$ is the incidence matrix $N$ of
a symmetric 2-$(36,15,6)$ design
(see  \cite[Theorem 2.4.12, page 99]{Ton88} or \cite[Sec. 3]{Ton98}). 
Replacing all  zeros in $N$
with $-1$'s gives a regular Hadamard matrix of order 36.
Using the clique finding algorithm Cliquer \cite{cliquer},
a quick computer search shows that the graph $\Gamma$ contains exactly 272
cliques of size 36, or in other words, there are 272 collections of 36 codewords
from $T$ that form a Hadamard matrix of order 36.
Further analysis with Magma shows that all 272 Hadamard matrices are equivalent
to a matrix $H'$ with a full monomial automorphism  group of order 72 that acts
transitively on the set of size 72 being  the union of the  rows of  $H'$ and the rows of 
$-H'$.

The incidence matrix of the symmetric 2-$(36,15,6)$ design $D'$ obtained  
by replacing all $-1$-entries of $H'$  with zeros  is 
listed in the Appendix. The design $D'$ has a trivial full automorphism group 
of order 1.

(ii) The 3-rank of the incidence matrix  of $D'$ is 18, and its row span
over $GF(3)$ is a ternary $[36,18,12]$ code equivalent to the
Pless symmetry code.

Parts (iii) was verified by computer using Magma.
The full automorphism group of $L(17)$ partitions the set of the
888 codewords of weight 36 into two orbits, of length 72 and 816
respectively, the orbit of length 72 comprised of the rows of $H$ (\ref{eq6})
and $-H$. Thus, the full automorphism group of the code $L(17)$
coincides with the full automorphism group of $H$, and is of 
order\footnote{This is  the order of the
Paley-Hadamard matrix of Type II for $q=17$ \cite{LS}.}
 $19584 =2^7 \cdot 3^2 \cdot 17.$
$\Box$

\begin{note}
{\rm
Up to equivalence, there are exactly 11 Hadamard matrices of order 36
with automorphism groups of order divisible by 17 (Tonchev \cite{Ton86}).
Each of these matrices spans a ternary self-dual code of length 36, but
only the symmetry code $C(17)$ spanned  the Paley-Hadamard matrix of type II is extremal, that is, has minimum distance 12, and supports 5-designs.
A stronger characterization of the Pless symmetry code $C(17)$ was proved
by Huffman \cite{Huf}, namely that  up to equivalence, $C(17)$ is the
only extremal ternary self-dual code of length 36 that admits a monomial automorphism of order $17$.
}
\end{note}
\begin{note}
{\rm
Hadamard matrices and designs are used
for the construction of self-orthogonal and self-dual codes 
over other finite fields. A classical example is the extended binary Golay code
generated by a  bordered incidence matrix of a symmetric Hadamard
2-$(23,11,5)$ design associated with a Paley-Hadamard matrix of type I.
Hadamard matrices of order 28 with an automorphism of order 7 \cite{Ton85}
were used by Pless and Tonchev \cite{PT} for the classification of self-orthogonal
codes over $GF(7)$. 
The Paley-Hadamard matrix of type II of order 28 is the only Hadamard matrix
of this order that admits an automorphism of order 13 and yields an extremal
binary self-dual code of length 56
\cite{Ton83}, \cite{Ton89}. More extremal binary self-dual codes derived from
Hadamard matrices of order 28 were found in \cite{BusT}.
}
\end{note}

\section{Hadamard matrices  and designs arising from ternary QR codes}
\label{s3}

The symmetry codes $C(11)$, $C(23)$, and $C(29)$
have siblings with the same parameters and weight distribution,
being ternary extended quadratic-residue codes
 that support 5-designs by  the Assmus-Mattson theorem.
 If $q\equiv 3 \pmod 4$
 is a prime power, a  quadratic residue (QR) code of length $q$ is  a code 
 spanned by the
(0,1)-incidence matrix $A$ of a symmetric Hadamard 2-$(q, (q-1)/2, (q-3)/4)$ design
 obtained from the Paley-Hadamard matrix of type I, and its extended
 code is spanned by a matrix obtained by adding one all-one column to $A$.
 If,  in addition,  $q\equiv -1 \pmod 3$, that is, $q$ is of the form $q=12s+11$
for some integer $s\ge 0$,  
 the ternary extended QR code is self-dual. 
\begin{thm}
\label{thm}
\label{t4}
Let $q=12s+11$ be a prime power, and let $QR_q$ be the ternary extended QR code of length $q+1$.\\
(i)  $QR_q$  contains a Paley-Hadamard matrix
of type I having as rows codewords of weight $q+1$.\\
(ii)  $QR_q$ contains a set of $2q$ (0,1)-codewords of weight $(q+1)/2$
that form the incidence matrix of a Hadamard 3-$(q+1,(q+1)/2, (q-3)/4)$ design
associated with the Paley-Hadamard matrix of type I of order $q+1$.\\
(iii) If $q=11$, 23 or 47, $QR_q$ contains exactly $2q$ (0,1)-codewords of weight 
$(q+1)/2$, and the permutation automorphism group of the code coincides
with the full automorphism group of the Hadamard 3-$(q+1,(q+1)/2, (q-3)/4)$ design from part (ii).
\end{thm}

{\bf Proof}. (i)  The statement (i) is implicit  in \cite{AM}.
 The column sum of the (0,1)-incidence $A$ of the Hadamard 
 2-$(q, (q-1)/2, (q-3)/4)$ design
 obtained from the Paley-Hadamard matrix of type I is
 \[ (q-1)/2 =(12s+10)/2 \equiv -1 \pmod 3, \]
 and the sum of all $q$ entries of the all-one column is $q\equiv -1 \pmod 3$,
 hence the sum over $GF(3)$ of all rows  of the $q\times (q+1)$-matrix $B$ obtained by bordering $A$ with one all-one column,  is equal to the constant vector  $\bar{2}$. Since  $QR_q$  is the row span of $B$,
 the constant vectors $\bar{1}$ and $\bar{2}$ belong to the code.
Let $E$ be the $(q+1)\times(q+1)$ matrix obtained from $B$ by adding one extra all-one row. The matrix $H_{q+1}=2J-E$, where $J$ is the 
$(q+1)\times (q+1)$ all-one matrix, is a Paley-Hadamard matrix of type I.
Every row of  $H_{q+1}$ is the difference of the codeword $\bar{2}$
and a row of $B$, hence the rows of $H_{q+1}$ belong to the code $QR_q$.

(ii)  Adding the codeword $\bar{2}$ to every row of
\begin{equation}
\label{eq}
 \left(\begin{array}{rr}
H_{q+1}\\
-H_{q+1}
\end{array}\right)
\end{equation}
gives a $(2q+2)\times (q+1)$ matrix with one all-zero row, one all-one row,
and $2q$ (0,1)-rows of weight $(q+1)/2$ that form the incidence matrix of a 
Hadamard 3-$(q+1,(q+1)/2, (q-3)/4)$ design associated with $H_{q+1}$.

(iii) The proof is similar to that of Corollary  \ref{cor1}.
$\Box$

\begin{note}
{\rm
The number of codewords of weight 60 in $QR_{59}$ is 41184. It seems likely
that there may be Hadamard matrices of order 60 formed by codewords of
weight 60 that are not equivalent to the Paley-Hadamard matrux of type I
from Theorem \ref{t4}.
}
\end{note}


\section{Concluding remarks}

The extended ternary Golay code of length 12, the 
Pless symmetry codes
$C(q)$ ($q=11, 17, 23$ and 29),
of lengths 24, 36, 48 and 60,
the extended ternary QR codes
of lengths 24, 48 and 60,
and an extremal code of length 60 discovered by Nebe and Villard
\cite{NV} as an analogue of the Pless symmetry code $C(29)$,
are the only known extremal ternary self-dual codes of length divisible by 12
that support 5-designs. It is known that the symmetry code 
of length 84 ($q=41$), as well as the extended QR code of this length are not extremal. Extremal ternary self-dual codes of length $n$ divisible 
by 12 do not exist  for $n=72, 96, 120$, and all $n\ge 144$, because then
the extremal Hamming weight enumerator contains a negative coefficient
\cite{RS}.  

All ternary self-dual codes of length 24 have been classified
up to equivalence (Harada and Munemasa \cite{HM}), and
the symmetry code $C(11)$ and the extended QR code 
are the only extremal  codes of this length.
Nine of the self-dual ternary codes of length 24 are spanned by
Hadamard matrices of order 24 \cite{LTP}, \cite{LPS},
but only  two codes, $QR_{23}$ and $C(11)$, that are spanned by the
Paley-Hadamard matrices of type I and II respectively, are extremal.

It is an interesting open question whether the Pless symmetry codes
of length 36, 48, and 60,  the extended QR codes of lengths 48 and 60,
and the extremal code of length 60 found by Nebe and Villard
\cite{NV}
are the only extremal self-dual codes of these lengths.
The results from Section \ref{s2} show that the symmetry code
of length 36 can be obtained  from a Hadamard matrix that is not a
Paley-Hadamard matrix of type II, and a natural question that arises is
whether any other extremal codes of length 36, 48, or 60 can be obtained
form Hadamard matrices that are not of Paley type.

The extremal ternary self-dual codes of lengths $n\ge 36$ have not been classified up to equivalence.  A partial classification
of such codes of length $n\le 40$ admitting  automorphisms of prime order 
$p\ge 5$ was given by Huffman \cite{Huf}.
 In addition, it was proved by Nebe \cite{Nebe}  that,  up to equivalence, the only extremal ternary self-dual codes
of length 48 that admit an automorphism of a
prime order $p\ge 5$, are the  Pless symmetry code and the 
extended QR code.

\section{Appendix}
\[
\begin{array}{c}
 1 1 0 1 1 0 0 0 1 1 0 1 0 0 1 1 0 0 0 0 1 0 1 0 1 0 0 0 0 0 1 1 0 0 0 1\\
 0 0 0 0 1 0 1 0 1 0 0 0 0 0 1 1 0 0 1 0 0 1 1 1 0 0 1 0 1 1 0 0 1 1 0 1\\
 0 1 1 0 0 0 1 0 0 0 1 0 0 0 1 1 0 0 0 1 1 0 0 1 1 1 0 0 1 1 0 1 0 0 1 0\\
 0 1 0 0 1 0 0 0 0 0 0 0 0 1 0 0 1 1 1 1 0 1 0 0 1 1 0 0 1 0 1 1 0 1 0 1\\
 0 0 1 0 1 1 1 1 0 1 0 0 1 0 0 0 1 0 0 1 1 1 1 1 1 0 0 0 0 0 0 0 0 0 0 1\\
 0 1 1 1 1 1 0 0 0 1 1 0 0 1 1 0 0 0 0 0 0 1 1 0 0 0 0 1 1 0 0 0 0 1 1 0\\
 0 0 0 0 0 0 0 0 1 1 1 1 1 1 0 0 0 0 1 1 0 1 1 0 1 0 0 0 0 1 0 1 1 0 1 0\\
 0 0 0 0 1 0 0 1 0 1 1 1 1 0 1 0 0 1 0 0 0 0 0 0 1 1 1 1 1 1 0 0 0 0 0 1\\
 0 1 0 0 0 1 0 0 0 1 0 1 1 1 1 1 1 0 0 1 0 0 0 1 0 0 1 0 0 1 1 0 0 1 0 0\\
 1 0 0 0 1 1 0 0 0 1 1 0 0 0 0 1 1 0 1 1 1 0 0 0 0 0 1 1 1 0 0 1 1 0 0 0\\
 0 0 0 1 1 1 1 1 1 0 0 0 0 0 0 0 0 0 0 1 0 0 0 0 1 0 1 1 0 1 1 1 0 1 1 0\\
 0 0 0 1 1 1 1 0 1 0 1 1 1 1 0 0 0 0 1 0 1 0 0 1 0 1 0 0 1 0 1 0 0 0 0 0\\
 0 1 0 1 0 0 1 1 0 0 1 0 1 0 1 1 1 0 1 0 0 1 0 0 1 0 0 1 0 0 1 0 1 0 0 0\\
 1 0 0 1 0 1 0 0 0 0 0 0 1 0 1 0 0 1 0 1 1 1 0 0 0 0 0 0 1 1 1 0 1 0 1 1\\
 1 0 1 0 0 1 0 0 1 0 0 1 0 0 1 0 1 0 1 0 1 1 0 0 1 1 0 1 0 1 0 0 0 1 0 0\\
 0 0 1 0 0 1 1 0 0 0 1 1 0 0 1 0 1 1 1 0 0 0 1 0 0 0 1 0 0 0 1 1 0 0 1 1\\
 0 1 0 0 0 1 1 1 1 1 0 0 0 1 1 0 0 1 0 0 1 1 0 0 0 1 1 0 0 0 0 1 1 0 0 0\\
 1 0 0 1 0 0 0 1 0 0 1 0 0 1 1 0 0 1 1 1 1 0 1 1 1 0 1 0 0 0 0 0 0 1 0 0\\
 1 0 0 0 1 0 1 0 1 0 0 0 1 1 1 1 1 1 0 1 0 0 1 0 0 1 0 1 0 0 0 0 0 0 1 0\\
 0 0 0 1 0 0 0 1 0 0 0 1 1 0 0 1 1 0 0 0 1 1 1 0 0 1 1 0 1 0 0 1 0 1 1 0\\
 1 1 1 1 1 0 0 0 1 0 1 0 1 0 0 0 1 1 0 0 0 1 0 1 0 0 1 0 0 1 0 1 0 0 0 0\\
 0 0 1 0 1 0 0 0 0 0 0 1 0 1 0 1 0 1 0 0 1 1 0 1 1 0 1 1 0 0 1 0 1 0 1 0\\
 1 1 1 0 0 1 0 1 1 0 0 0 1 1 0 1 0 0 1 0 0 0 0 0 1 0 1 0 1 0 0 0 0 0 1 1\\
 1 1 0 1 0 0 1 0 0 1 0 1 0 0 0 0 0 0 1 1 0 1 0 1 0 1 1 1 0 0 0 0 0 0 1 1\\
 0 0 0 1 0 1 0 0 1 1 1 0 0 0 0 1 1 1 0 0 0 0 0 1 1 1 0 0 0 0 0 0 1 1 1 1\\
 1 0 1 0 0 0 1 1 1 1 1 1 0 0 0 1 0 1 0 1 0 1 0 0 0 0 0 0 1 0 1 0 0 1 0 0\\
 0 1 0 0 0 0 0 1 1 1 0 0 0 0 0 0 1 1 1 0 1 0 1 1 0 0 0 1 1 1 1 0 0 0 1 0\\
 0 0 1 1 0 0 1 0 0 1 0 0 1 1 0 1 0 1 1 0 1 0 0 0 0 0 0 1 0 1 0 1 0 1 0 1\\
 0 0 1 1 0 0 0 1 1 0 0 1 0 1 1 0 1 0 0 1 0 0 0 1 0 0 0 1 1 0 0 1 1 0 0 1\\
 0 1 1 0 0 0 0 0 1 0 1 0 1 0 0 0 0 0 0 1 1 0 1 0 0 1 1 1 0 0 1 0 1 1 0 1\\
 1 0 0 0 0 1 0 1 0 0 1 0 0 1 0 1 0 0 0 0 0 1 1 1 0 1 0 1 0 1 1 1 0 0 0 1\\
 1 0 1 0 1 0 0 1 0 1 0 0 1 0 1 0 0 0 1 0 0 0 0 1 0 1 0 0 0 0 1 1 1 1 1 0\\
 0 1 1 1 1 1 0 1 0 0 0 1 0 0 0 1 0 1 1 1 0 0 1 0 0 1 0 0 0 1 0 0 1 0 0 0\\
 1 0 1 1 0 0 1 0 0 1 0 0 0 1 0 0 1 0 0 0 0 0 1 0 1 1 1 0 1 1 1 0 1 0 0 0\\
 1 1 0 0 1 0 1 1 0 0 1 1 0 1 0 0 1 0 0 0 1 0 0 0 0 0 0 0 0 1 0 0 1 1 1 1\\
 1 1 0 0 0 1 1 0 0 0 0 1 1 0 0 0 0 1 0 0 0 0 1 1 1 0 0 1 1 0 0 1 1 1 0 0\\
 \end{array}
\]
A $2$-$(36,15,6)$ design associated with the Pless symmetry code
of length 36

\end{document}